\documentclass[12pt]{article}

\usepackage{amsmath,amsthm,amsfonts,latexsym,amscd}

\usepackage{eucal}

\usepackage{epsfig,graphicx}

\swapnumbers

\theoremstyle{plain}
              \newtheorem{theorem}{Theorem}[section]
              \newtheorem{lemma}[theorem]{Lemma}
              \newtheorem{corollary}[theorem]{Corollary}
              \newtheorem{proposition}[theorem]{Proposition}

\theoremstyle{definition}

              \newtheorem{notation}[theorem]{Notation}
              
              \newtheorem{remark}[theorem]{Remark}

\numberwithin{equation}{section}

\addtocounter{section}{-1}

\newcommand{\script}[1]{{\mathcal{#1}}}

\newcommand{\gs}{\sigma}
\newcommand{\ga}{\alpha}
\newcommand{\gb}{\beta}
\newcommand{\gl}{\lambda}

\newcommand{\eps}{\epsilon}
\newcommand{\ta}{\theta}

\newcommand{\xv}{\mathbf x}
\newcommand{\yv}{{\mathbf y}}
\newcommand{\tv}{\mathbf t}

\newcommand{\uv}{\mathbf u}

\newcommand{\zv}{\mathbf z}

\newcommand{\qtspp}[1]{q_{\tv_{#1},s_{#1}}^+}
\newcommand{\qtsmm}[1]{q_{\tv_{#1},s_{#1}}^-}

\newcommand{\bt}[1]{b_{\tv_{#1}}}
\newcommand{\ptsp}{p_{\tv,s}^+}

\newcommand{\ptsm}{p_{\tv,s}^-}

\newcommand{\ptspm}{p_{\tv,s}^{\pm}}

\newcommand{\rn}[1]{$\mathbb R^{#1}$}

\newcommand{\conv}{\operatorname{conv}}

\newcommand{\lan}{\langle}
\newcommand{\ran}{\rangle}

\newcommand{\pw}[1]{\partial W_{#1}}

\DeclareMathOperator{\seg}{Seg}

\DeclareMathOperator{\range}{Range}
\DeclareMathOperator{\sg}{Seg}

\begin{document}
\title{The Spectral Scale and the Numerical Range}
\author{Charles A. Akemann and Joel Anderson
\thanks{The second author was partially supported by the National
Science Foundation during the period of research that resulted in this
paper.}}

%\noindent
%\address{Department of Mathematics\\
%University of California\\
%Santa Barbara, CA 93106}

%\email{\tt akemann\@math.ucsb.edu}

%\noindent
%\address{Department of Mathematics\\
%University Park, PA 16802}

%\email{\tt anderson\@math.psu.edu}

\maketitle
\begin{abstract}
Suppose that $c$ is an operator on a Hilbert Space $H$ such that the
von Neumann algebra $N$ generated by  $c$ is finite.  Let $\tau$ be a faithful
normal tracial state on $N$ and set $b_1 = (c+c^*)/2$ and $b_2 = (c-c^*)/2i$.
Also write $B$ for the spectral scale of $\{b_1, b_2\}$ relative to $\tau$.
In previous work by the present authors, some joint with Nik Weaver,  $B$ has
been shown to contain considerable spectral information about the operator
$c$. In this paper we expand that information base by  showing that the
numerical range of $c$ is encoded in $B$ also.

We begin by proving that the $k$-numerical range of an arbitrary operator $d$
in $B(H)$ coincides with the numerical range of $d$ when the von Neumann
algebra generated by $d$  contains no finite rank operators.  Thus, the
$k$-numerical range is not useful for most operators considered here.

We next show that the boundary of the numerical range of $c$ is exactly the set
of radial complex slopes on $B$ at the origin.  Further, we show that points on
this boundary that lie in the numerical range are visible as line segments in
the boundary of
$B$. Also, line segments on  the boundary which lie in the numerical range show
up as faces of dimension two in the boundary of $B$.  Finally, when $N$ is
abelian, we prove that the point spectrum of $c$ appears as complex slopes of
1-dimensional faces of $B$.

AMS Subject Classification Numbers 47A12, 47C15

\end{abstract}

\section{Introduction and Notation}
\medskip
We shall develop notation here that will be used throughout the
rest of the paper. Suppose that $c$ is an operator on a  separable
infinite dimensional Hilbert Space $H$ such that the von Neumann algebra $N$
generated by  $c$ is finite and contains the identity operator of $H$.  Let
$\tau$ be a faithful normal tracial state on $N$ and set $b_1 = (c+c^*)/2$ and
$b_2 = (c-c^*)/2i$.  Also write $B = B(b_1,b_2)$ for the spectral scale of
$\{b_1, b_2\}$ relative to
$\tau$.   Recall from \cite{AAW} that the spectral scale is defined as follows.
We define a map  $\Psi$ on $N$ by  the formula
\[
\Psi(a) = (\tau(a),\tau(ab_1), \tau(ab_2))
\]
and write $B =\Psi(N_1^+)$, where $N_1^+ = \{a\in N: 0 \le a \le 1\}$.
For the purpose of the present paper it is more natural to identify the second
and third real coordinates of $\Psi(a)$ with a point in $\mathbb C$ and  to
view the range of
$\Psi$ as lying in
$\mathbb R
\times
\mathbb C$ i.e.
\[
\Psi(a) = (\tau(a),\tau(ac))).
\]
           Since $\tau$ is normal, $B = B(c)$ is a compact, convex subset
of $\mathbb
R\times \mathbb C$, which we call the {\bf spectral scale} of $c$ {\bf relative
to $\tau$}. (Since $\tau$ is fixed throughout, we will suppress the
dependency on $\tau$ in the sequel.) This situation arises whenever
$N$ is finite dimensional.  That case is analyzed fully in \cite{numrange}.
The case of infinite dimensional $N$ arises  most naturally when $c$ is an
element in a factor of type $II_1$ with the unique normal trace
$\tau$.

In  \cite{numrange} we showed that if $c$ acts on finite dimensional Hilbert
space, then the spectral scale contains a canonical affine image of each of
the k-numerical ranges of $c$.  In the present paper  we restrict our attention
to infinite dimensions and show how to derive information about both the
spectrum  and the numerical range of the  operator $c$ from $B$.

A key notion is that of complex slope.  This is defined naturally in $\mathbb
R \times \mathbb C$ for the line segment joining two points
$(x_1,z_1),(x_2,z_2)$
to be
\[
\frac{z_2-z_1}{x_2-x_1}.
\]
As shown in \cite{AAW} and \cite{Geom II} the geometry of the spectral scale
reflects spectral properties of real linear combinations of $\{b_1, b_2\}$.
This information is contained in slopes of 2-dimensional projections of $B$, so
it is natural to look for numerical range information is the same way, i.e. as
slopes.  Just as  in \cite{AAW} the notion of slope needs to be extended beyond
segments to the slope  of a curve, i.e. a derivative.  In the present paper we
define the radial complex slopes of $B$ as a certain of complex directional
derivative in Section 2, and  prove that the set of such slopes is exactly
the boundary of the numerical range of  $c$.  As part of our analysis we show
that the numerical range of $c$ is exactly the set of complex slopes of line
segments in $B$ that are anchored at the origin.  In addition if $\gl$ lies on
the boundary of the numerical range, then the corresponding line segment lies
on the boundary of $B$.  Further, a line segment in the numerical range that
lies in its boundary corresponds to a 2-dimensional face of $B$ that contains
the origin.  Finally, when $N$ is abelian,  we show that the point spectrum of
$c$ appears as complex slopes of 1-dimensional faces of $B$.  Before we relate
$B$ to the numerical range of
$c$, we first prove some results about the numerical range and k-numerical
range of more general operators in
$B(H)$.

\section{The Numerical Range and the k-Numerical range}
\medskip

{\bf Notation for this section}: Fix an arbitrary element $d
\in B(H)$, write
$d_1=(d+d^*)/2 , d_2=(d-d^*)/2i$, let $A$ denote the C*-algebra generated by
$d$ and the identity $1$ of $B(H)$ and let $M$ denote the von Neumann algebra
generated by $d$ and
$1$.  The {\bf numerical range}
of
$d$ is by definition
\[
W(d) = \{\lan d\xv,\xv\ran: \xv\in H \text{ and } \|\xv\| =1\}.
\]

We first recall some known facts about the numerical range which we state as a
theorem for easy reference.

\begin{theorem} The following statements hold.
\begin{enumerate}

\item $W(d)$ is a bounded convex subset of $\mathbb C$, which is not
necessarily closed.

\item If $d$ is a diagonal operator in $B(H)$ with eigenvalues
$\gl_1,\gl_2,\dots$, then
\[
W(d) =
\left\{ \sum_{n=1}^\infty \gl_n t_n:  0 \le t_i \le 1 \text{
and }\sum_{n=1}^\infty t_n = 1\right\}.
\]
In other words, $W(d)$ is the (infinite) convex hull of the eigenvalues of $d$.
\item If $d = d^*$ and
\[
\gb^- = \inf\{\gb\in \gs(d)\} \text{ and } \gb^+ = \sup\{\gb\in\gs(d)\},
\]
then $(\gb^-,\gb^+) \subset W(d) \subset[\gb^-,\gb^+]$.

\item  If $W(d)$ is a line segment, $\gl$ is  an end point of $W(d)$ and
$\xv$ is a unit vector such that $\gl = \lan d\xv,\xv\ran$, then $d\xv =
\gl\xv$ and $\gl$ is a reducing eigenvalue for $d$.

\end{enumerate}

\end{theorem}
\begin{proof}  The first part of assertion $(1)$ is known as the {\em
Toeplitz--Hausdorff Theorem} (see \cite{Toep} and \cite{Haus}).  To
see that $W(d)$ is not necessarily closed, consider the infinite diagonal
matrix $b$  whose eigenvalues are $1,1/2,\dots, 1/n,\dots$.  Applying
part
$(2)$ of the Theorem, we get  $W(b) = (0,1]$.

For $(2)$, if $\xv = (x_1,x_2,\dots)$ is a unit vector in $H$, then we have
\[
\lan d\xv,\xv\ran = \sum_{n=1}^\infty \gl_n|x_n|^2.
\]
Putting $t_n = |x_n|^2$, we get that $W(d)$ is contained in the right hand side
of the formula in (2).  On the other hand, if $\sum_{n=1}^\infty
\gl_n t_n$ is as
on the right hand side of (2) and we set $\xv =(\sqrt t_1, ...)$, then $\lan
d\xv,\xv\ran =
\sum_{n=1}^\infty \gl_n t_n$ and so the reverse inclusion also holds.

Assertion (3) follows immediately from $(2)$ when $d$ is diagonal.  The general
result is obtained from this fact and spectral theory.

Translating, rotating and scaling if necessary, we may assume
that $(0,1) \subset W(d) \subset [0,1]$ and $\lan d\xv,\xv\ran = 1$.  In this
case $d= d^*$ and $\|d\| = 1$. We have
\[
\|(d-1)\xv\|^2 = \|d\xv\|^2 - \lan d\xv,\xv\ran - \lan \xv,d\xv\ran + 1 \le
\|d\| - 1 = 0.
\]
\end{proof}
\medskip

The $k$--numerical range was introduced by Berger in his thesis
\cite{Berg}.  It is defined by the formula
\[
W_k(d) = \left\{\frac{1}{k}\sum_{i=1}^k \lan
d\xv_i,\xv_i\ran:\text{ the $\xv_i$'s are orthonormal}\right\},\quad 1
\le k \le n.
\]
Since we have $W_1(d) = W(d)$,  this notion is a generalization of the
standard numerical range.  Berger showed in \cite{Berg} that the $k$-numerical
range is convex.  (See \cite[Problem 167]{Halmos} for a proof of this fact).

\medskip
\begin{theorem} Suppose $\iota$ denotes the canonical (unbounded)
trace  on $B(H)$
and $k$ is a positive integer.  If
\[
K_k = \{a \in B(H), 0 \le a \le 1, \iota(a) = k\}
\]
and
\[
V_k =\{\iota(ad)/k : a\in K_k\},
\]
then
\[
W_k(d) \subset V_k \subset \overline {W_k(d)}.
\]
\end{theorem}
\begin{proof}  Fix an element $\gl \in W_k(d)$ so that $\gl$ has the form
\[
\gl = \frac{1}{k}\sum_1^k\lan d\xv_i,\xv_i\ran,
\]
where the vectors $\xv_1,\dots,\xv_k$ are orthonormal. If we write
$p$ for the projection onto the span of $\{\xv_1, ..., \xv_k\}$, then $p$ has
rank $k$ so that $\iota(p) = k$.  Hence, $p \in K_k$ and we have $\iota(pd)/k
=(1/k)\sum_1^k\lan d\xv_i,\xv_i\ran = \gl$.  Thus, $W_k(d) \subset V_k$.  Note
that this calculation also shows that if $p$ is a projection of rank $k$, then
$\iota(pd)/k \in W_k(d)$.

Now observe that the finite rank operators in $K_k$ are dense in $K_k$ with the
trace  norm, and $\iota$ is continuous for that norm.  Thus, to establish the
second inclusion,  it suffices to take a finite rank element $a\in
K_k$ and show
that $\iota(ad)/k  \in W_k(d)$.  So suppose that $a$ is such an element and let
$H_0$ denote the range of $a$ so that we may view $a$ as acting on
$B(H_0)$.  We
have then that $a$ is a convex combination of projections of rank $k$
in $B(H_0)$
by \cite[Corollary 1.2]{numrange}. Since
$W_k(d)$  is convex, the second inclusion follows from the last sentence of the
previous paragraph.
\end{proof}
\medskip

The numerical range and the k-numerical range of $d$ are dependent on the
way $d$ acts on a Hilbert space, i.e. they are operator theory concepts
rather than operator algebra concepts.   In order to get a set which is tied
to the C*-algebra $A = \text{C}^*(d,1)$, we define the abstract
numerical range
of
$d$, which we denote by $W_{ab}(d)$.  This set is defined  by the formula
\[
W_{ab}(d) = \{f(d):f\text{ is a state on $A$} \}.
\]
The fact that $W_{ab}(d) =\overline{W(d)}$ is known.  We
include a proof for completeness.
\begin{proposition}
\[
W(d)^\circ = W_{ab}(d)^\circ \text{ and } \overline{W(d)} =
W_{ab}(d).
\]
\end{proposition}
\begin{proof}
It suffices to show that the second assertion holds since the first
assertion then follows from the convexity of  $W(d)$ and
$W_{ab}(d)$ by \cite[2.3.8]{Web}. Since the set $S$ of states on $A$  is convex
and compact in the weak* topology, $W_{ab}(d)$ is  convex and closed
and since
\[
W(d) = \{f(d): f \text{ is a vector state of $A$}\}
\]
we get $W(d)\subset W_{ab}(d)$.  The proof that $\overline{W(d)} =
W_{ab}(d)$ is
completed by noting that the convex hull of the
vector states of $A$ is weak$^*$--dense in $S$ by
\cite[3.4.1]{Dix}.
\end{proof}

In most of the rest of this section and the next, we shall be
interested in studying points in the relative boundary of $W_{ab}(d)$. In
particular we will be interested  when  points in this set lie in $W(d)$.
This study is facilitated by introducing some new notation and making a
normalization.

\begin{notation}
If $S$ is a convex subset of $\mathbb C$, we write $\partial S$ for
its relative
boundary.   Fix $\gl \in \pw{ab}(d)$ and let $F_{ab}$ denote the (proper) face
of maximal dimension in $W_{ab}(d)$ that contains $\gl$.  Translating and
rotating if necessary, we may assume that  $\gl = 0$, $W_{ab}(d)$ lies in the
right half plane and $F_{ab}$ lies on the imaginary axis.  Now write $F_W =
F_{ab}\cap W(d)$.  Observe that $F_W$ may be empty, but if $F_W \ne \emptyset$,
then it is a face in $W(d)$.    Finally let $r$ denote the projection onto
$\script N(d_1)$ (the null space of $d_1$), write $p$ for the projection onto
$\script N(rd_2r)$ and set $q  = rp$.
   \end{notation}

\begin{lemma} With the notation and normalization introduced in 1.4 above, the
following statements hold.

\begin{enumerate}

\item $F_W \ne \emptyset$ if and only if $r > 0$.

\item  $F_W = \{0\}$ if and only if  $0 < q =r$.

\item  $F_W$ is one dimensional and $0$ is an endpoint of $F_W$ if and only if
$rd_2r$ is semi-definite.  Further, if $F_W$ is one dimensional, $0$ is an
endpoint of $F_W$ and $0 \in W(d)$, then $0 < q < r$.

\item $F_W$ is one dimensional and $0$ is in the interior of $F_W$ if
and only if
$(rd_2r)^\pm \ne 0$.
\end{enumerate}
\end{lemma}

\begin{proof}  Since $F_W$ is the intersection of $W(d)$ with the imaginary
axis, we have
\[
F_W = \{\lan d\xv,\xv\ran:\|x\| = 1\text{ and } \lan d_1 \xv,\xv\ran = 0\}.
\]
Thus, $F_W\ne \emptyset$ if and only if there is a unit vector $\xv$ such that
$\lan d_1 \xv,\xv\ran = 0$.  Since $W_{ab}(d)$ lies in the right half plane,
$d_1 \ge 0$ and so $\lan d_1 \xv,\xv\ran = 0$ if and only if $d_1\xv =0$ by
part (4) of Theorem 1.1.  In other words $F_W \ne \emptyset$ if and only if
$r > 0$.   Thus, (1) holds.

Next, the calculation above also shows that $F_W = \{0\}$ if and only if for
each unit vector $\xv$ in the range of $r$ we have $\lan d_2\xv,\xv\ran = 0$
and this occurs if and only if $rd_2r = 0$, which in turn occurs if and only
if $0 < q = r$. Thus, (2) is true.

Now suppose $F_W$  is one dimensional.  In this case $0$ is an endpoint if and
only if $F_W$ lies on the positive or negative imaginary axis and this occurs
if and only if $rd_2r$ is semi-definite.  Thus the first assertion in (3)
holds.  For the second assertion, suppose that $0$ is an endpoint of $F_W$ and
$0\in W(d)$. In this case  there is a unit vector $\xv$ such that  $\lan
d_1\xv,\xv\ran = \lan d_2\xv,\xv\ran = 0$ and since $d_1 \ge 0$ and $rd_2r$ is
semi-definite, we must have $d_1\xv = d_2\xv = 0$ and therefore $0 < q \le r$.
Since $F_W \ne \{0\}$, $q< r$by part (2) of the Lemma and so (3) holds.
Finally, $0$ is in the interior of $F_W$ if and only if $rd_2r$ is not
semi-definite; i.e., if and only if $(rd_2r)^\pm \ne 0$. Thus (4) holds.

\end{proof}

The next result shows that for many interesting operators, the k-numerical
range is the same as the numerical range.  This suggests that the spectral
scale may be the more useful object for such operators because its geometry
displays the numerical range  (as we show in the next section) and much
more.

\begin{theorem} If each projection in $M=\{d,1\}^{\prime\prime}$ has infinite
rank in $B(H)$, then $W_k(d)=W(d)$ for each $k$.
\end{theorem}
\begin{proof}  First note that  $W_k(d) \subset W(d)$  by the definition of
$W_k(d)$ and the fact that $W(d)$ is convex.  In order to establish the reverse
inclusion,  we first show that the following
statements hold.

\smallskip

{\narrower\narrower

\noindent 1. If $\gl$ is in the relative boundary of $W(d)$ and $\gl$ is in
$W(d)$, then $\gl$ is in $W_k(d)$.

\noindent 2.  If $\gl$ is in the relative boundary of $W(d)$, but
$\gl$ is not in $W(d)$, then $\gl$ is in the relative boundary of $W_k(d)$, but
$\gl$ is not in $W_k(d)$.

}

\smallskip

Fix $\gl$ in the relative boundary of $W(d)$.   Normalizing as in 1.4 above,
we may assume that $\gl = 0$ and $W_{ab}(d)$ lies in the right half plane.
Also, we may define $F_W$,  $r$ and $q$ as in 1.4.

\medskip

\noindent {\bf Proof of statement 1:} Since $\gl = 0$ is in $W(d)$, we have
$F_W \ne \emptyset$ so that $r> 0$ by part (1) of Lemma 1.5.  If $F_W = \{0\}$
or $0$ is an endpoint of $W(c)$, then $q > 0$ by parts (2) and (3) of
Lemma 1.5.
Since $q$ has infinite rank by hypothesis, we may find an infinite orthonormal
set $\{\xv_n\}$ in the range of $q$ such that
\[
\lan d\xv_n,\xv_n\ran = 0 \quad n = 1,2,\dots.
\]
Thus $0 \in W_k(d)$ for all $k$ in this case.

Now suppose $F_W$ is one dimensional and  $0$ is in the interior of $F_W$.  In
this case, then we get $(rd_2r)^\pm \ne 0$ by part (4) of lemma 1.5.  Since the
range projections  of $(rd_2r)^\pm$ are infinite dimensional by hypothesis,
there exist mutually orthogonal  infinite orthonormal sets $\{\xv_j\}$ and
$\{\yv_j\}$ in $rH$ such that for each j, $\lan d_2\xv_j,\xv_j \ran > 0$ and
$\lan d_2\yv_j,\yv_j \ran < 0$.  Hence, for each j the subspace spanned by
$\xv_j$ and $\yv_j$ contains a unit vector $\zv_j$ such that $\lan
d_2\zv_j,\zv_j \ran = 0$ and the argument may now be completed as above to show
that $0 \in W_k(d)$ in this case.  Hence, statement 1 above holds.

\medskip

\noindent {\bf Proof of statement 2:}
Now suppose $\gl = 0$ is in the relative boundary of $W(d)$, but it is not in
$W(d)$ (so that 0 is not in $W_k(d)$). Also, we continue to assume  that no
point of $W(d)$ has negative real part so that $d_1 \ge 0$.  Since $0 \in
W_{ab}(d)$,  there is a state $g$ of $A = \text{C}^*(d,1)$ such that $g(d)=0 =
g(d_1) = g(d_2)$ by Proposition 1.3.  We can extend $g$ to a state of $M$,
denoted also by $g$. Write $p_n$ for the spectral projection of $d_1$
(computed in $M$) corresponding to the characteristic function of the interval
$(1/n,\infty)$ and set $r_n = 1-p_n$.  Since $g(d_1)=0$ and $d_1\ge 0$, it
follows that $g(p_n)=0$ by \cite[pages 304 and 305]{And}, so that $g(r_n)=1$,
for all n. Now consider $r_nd_2r_n$. Since $g(r_n)=1$, we get that
$g(r_nd_2r_n) = g(d_2) = 0 $ for all n. Thus $\overline{W(r_nd_2r_n)}$
contains 0 for all $n$.   We now consider two sub-cases.

$1^\circ$ Suppose that the null space of $d_1$ is 0 so that the $r_n$'s
decrease to 0 and fix $\epsilon > 0$.  In this case we may  choose a
unit vector
$\yv_1
\in r_1H$ such that
$|\lan r_{1}d_2r_{1}\yv_1,\yv_1 \ran| < \epsilon$ and set $j_1=1$.  Since
$\|r_n\yv_1\| \rightarrow 0$ as $n\rightarrow \infty$, we can find an index
$j_2$ so large that $r_{j_2} < r_{j_1}$ and if we set
\[
\xv_1 = \frac{1}{\|(r_{j_1}-r_{j_2})\yv_1\|}(r_{j_1}-r_{j_2})\yv_1,
\]
then
\[
|\lan r_{j_1}d_2r_{j_1}\xv_1,\xv_1 \ran| < \epsilon.
\]
Observe that if  $n \ge j_2$, then  $r_n \le r_{j_2}$ and  $r_n\xv_1 = 0$.

Now let us proceed by induction.  Assume that for some $m \ge 1$ we have
chosen orthogonal unit vectors $\{\xv_1, ..., \xv_m\}$ and  natural numbers
$1=j_1 < j_2 < ...< j_{m+1}$ such that for $i=1, ..., m$
\[
\xv_i=(r_{j_i}-r_{j_{i+1}})(\xv_i) \text{ and } |\lan d_2\xv_i,\xv_i \ran| <
\epsilon
\]
Since $\overline{W(r_{j_{m+1}}d_2r_{j_{m+1}})}$ contains 0, we can choose a
unit vector $\yv_{m+1} \in r_{j_{m+1}}H$ such that $|\lan
d_2\yv_{j_{m+1}},\yv_{j_{m+1}} \ran| < \epsilon$.  Arguing as above we may find
a unit vector $\xv_{j_{m+1}}$ and an index $j_{m+2}$  such that  $(r_{j_{m+1}}-
r_{j_{m+2}})\xv_{j_{m+1}}=\xv_{j_{m+1}}$ and $|\lan d_2
\xv_{m+1},\xv_{m+1}\ran| < \epsilon$.  Hence, the induction continues.

This process produces an infinite sequence $\{\xv_n\}$ of orthonormal vectors
such that $|\lan d_2\xv_n,\xv_n \ran| <\eps$  for each $n$.  Since $|\lan
d_1\xv_n,\xv_n\ran| \le 1/n$ by the definition of the $r_n$'s, we get
that $|\lan
d\xv_n,\xv_n\ran| < 2\epsilon$ for all large $n$. Thus by selecting $n$ large
enough, and relabeling, we get orthonormal vectors $\xv_1,\dots.
\xv_k$ such that
\[
|\lan d\xv_i,\xv_i\ran| \le 2\eps, \quad i = 1,\dots,k.
\]
Since $\epsilon$ was arbitrarily small, we get that $0 \in
\overline{W_k(d)}$ in
this sub-case.

$2^\circ$ Suppose that the null space of $d_1$ is nonzero and let $r$ denote
the projection onto this null space.  As above, we have $g(r) = 1$ and $g(rd)
= g(d_1) = 0$  Since the closure of $W(rd_2r)$ contains 0 and $r$ is infinite
dimensional,  we may select a select an orthonormal sequence $\xv_n$ of unit
vectors in the range of  $r$ such that
$\lan d\xv_n,
\xv_n\ran = \lan d_2\xv_n, \xv_n\ran < 1/n$.  Hence, as above, we get $0 \in
\overline W_k(d)$.  Hence in all cases $\gl$ is in the relative boundary of
$W_k(d)$.  Thus, statement 2 above holds.

Hence, the relative boundary of $W(d)$ is contained in the relative boundary of
$W_k(d)$. Since we also know that $W_k(d)$ is contained in $W(d)$, it follows
that $W(d)$ and $W_k(d)$  have the same relative boundaries.  Hence, by
\cite[2.3.8]{Web} these sets have the same relative interiors.  Finally, we get
that a point on this joint relative boundary lies in $W(d)$ if and only if it
lies in $W_k(d)$  by statement 1 above and the containment of $W_k(d)$ in
$W(d)$.  Hence these sets are equal.
\end{proof}

The next result is probably known, although we have not been able to find a
reference.  It is included here because it will be used in the proof of
Corollary 2.12 below.

\bigskip

\begin{theorem} If $d$ is unitary, then $\sigma(d)$ is exactly the
set of extreme points $z$ of $W_{ab}(d)$ such that $|z|=1$.
\end{theorem}
\begin{proof} We have that $\overline{W(d)} = W_{ab}(d)$ is contained in the
unit disk because $d$ is unitary.  Since $\gs(d)$ is contained in both the
unit circle and $W_{ab}(d)$ it follows that every point in $\gs(d)$ is an
extreme point of $W_{ab}(d)$.

For the reverse inclusion, recall that since $c$ is unitary, $\overline{W(d)} =
\conv(\gs(d))$ by \cite[Problem 171]{Halmos} and so the extreme
points of $\overline{W(d)}$ lie in $\gs(d)$.
\end{proof}

\begin{remark} The results in  Lemma 1.5 are related to some work of Gustafson
and Rao in \cite[\S 1.5]{Gust}, where they studied points in  $W(c)$ using
the set
\[
M_\gl = \{\xv: \lan c\xv,\xv \ran = \gl \|\xv\|^2\}, \quad \gl\in \mathbb C.
\]
Observe that $M_\gl \ne \{0\}$ if and only if $\gl \in W(c)$, but
that, in general this set is not a proper linear subspace.

Gustafson and Rao
showed in \cite[Theorems 1.5-1, 1.5-2 and 1.5-3]{Gust} that the following
statements hold.

{\smallskip \narrower\narrower
\noindent A.  If $\gl \in W(c)$, then $\gl$ is an extreme point of $W(c)$ if
and only if $M_\gl$ is a linear subspace.

\noindent B.  If $\gl$ is in $\pw{}(c)\cap W(c)$ and $\gl$ is not an extreme
point of $W(c)$ and $L$ is the the line of support for $\gl$, then
\[
\bigcup_{\ga\in L} M_\ga
\]
is a closed subspace of $H$.

\smallskip}

The fact that if $\gl$ is not an extreme point then $M_\gl$ is not a subspace,
follows from consideration of  $2\times 2$ matrices.  Indeed, the problem may
be reduced to showing that $M_0$ is not linear for the matrix
\[
\bmatrix
1&a\\
0&-1
\endbmatrix
\]
and this is easy to calculate.

The remaining assertions follow from Lemma 1.5.  In fact, if  we normalize as
in 1.4 so that $\gl = 0 \in \pw{}(c)\cap W(c)$ and the line of support  $L$ is
the imaginary axis, we use the other notation developed there, and suppose
that  $\gl = 0$ is an extreme point of $W(c)$,  then either $F_W = \{0\}$ or
$0$ is an endpoint of $F_W$.  If $F_W = \{0\}$, then $0 < q =r$ by part (2) of
Lemma 1.5 and we have $M_0 = \script N(b_1) = qH$.  If $0$ is an endpoint of
$F_W$, then since $0 $ is in $W(c)$ we have $0 < q < r$ by part (3) of Lemma
1.5 and we get $M_0 = \script N(q) $.   Thus $M_0$ is a subspace in both cases.

The second assertion follows from the fact that if  $F_W \ne \emptyset$
then
\[
\script N(b_1) = \bigcup_{\gl \in F_W} M_\gl,
\]
which is easily established using part (1) of Lemma 1.5 and part (4) of lemma
1.1.

\end{remark}

\bigskip
\section{Complex slopes}
\bigskip

In this section, we return to the notation developed in section 0. Thus, $N$
is the finite von Neumann algebra generated by the element $c$.  Our goal here
is to prove an analogue of Theorem 1.4 in \cite{numrange} in the infinite
dimensional case.  This will require some  preparation.

It is useful to begin by reviewing the situation when $c$ acts on a  Hilbert
space of dimension  $n$.  The key here is the fact that $N$ now has minimal
projections and the range of the trace on the projections in $N$ is finite.
It follows that the extreme points of $B$ have the form $(k/n,z)$, where $k$ is
an integer between $1$ and $n$. Thus if we define the $k^{\text{th}}$ {\em
isotrace slice} of $B$ by
\[
I_{k/n} = \{(k/n,z)\in B\},
\]
then the extreme points of $B$ lie in these isotraces slices and so
\[
B = \conv(0, I_{1/n}\cup \cdots \cup I_{k-1/n}\cup \Psi(1)).
\]
We showed in \cite[Theorem 1.4]{numrange} that  each $I_{k/n}$ is an affine
image of $W_k(c)$.  In particular, for $k = 1$,  the map
\[
\gl \mapsto \frac{1}{n}(1,\gl)
\]
is an affine isomorphism of $W(c)$ onto $I_{1/n}$.  Thus, if $(1/n,\gl/n)$
lies on the boundary of $I_{1/n}$ so that $\gl\in W(c)$, then the line
segment joining this point and the origin lies on the boundary of $B$.  The
complex slope of this map is
\[
\frac{\gl/n - 0}{1/n-0} = \gl.
\]
Further, if $F$ is a face of dimension one on the boundary of $W(c)$, then
$I_{1/n}$ has a corresponding one dimensional face on its boundary and the
line segments from points on this face to the origin form a face of dimension
two in the boundary of $B$.  We shall show below that the precise analogs of
these facts hold in infinite dimensions.

On the other hand the identification of $W(c)$ with $I_{1/n}$ does not
carry over as nicely.  This is a result of the fact that in infinite
dimensions we may  no longer have minimal projections.  In fact if we now
define  an {\bf isotrace slice of $B$} to be  a set of the form
\[
I_t = \{\xv = (t,z) \in B\},
\]
for $0 < t < 1$, then  the extreme points of $B$ may lie in a continuum of
isotrace slices.    Thus, it is impossible in
this case to identify the numerical range as a multiple of an isotrace slice
as is the case in finite dimensions. Nevertheless as we shall show below in
Theorem 2.4, the isotrace slices do determine the abstract numerical range  in
a way that generalizes the finite dimensional case.

To see how this identification arises, let us return for a moment to the
finite dimensional case.  As noted above, the portion of $B$ that lies between
the planes $x=0$ and $x = 1/n$ is the convex hull of 0 and $I_{1/n}$.  Hence,
if $0 < t < 1/n$, then $tI_{1/n} = I_{t/n}$.  Equivalently, we have
\[
I_{1/n} = \frac{I_{t/n}}{t}.
\]
Since
\[
W(c) = nI_{1/n} = \frac{I_{1/n}}{1/n}
\]
we get
\[
\frac{I_{t/n}}{t/n} = \frac{I_{1/n}}{1/n} = W(c), \qquad 0 < t \le \frac{1}{n}.
\]

Thus, in the general infinite dimensional case, it is natural to
consider the  map $\Delta :B(c)\setminus \{0\} \rightarrow \mathbb C$ defined
by
\[
\Delta(x_0, re^{i\theta}) = (r/x_0)e^{i\theta}.
\]
and to view sets of the form $\Delta_t = \Delta(I_t),\quad 0 < t < 1$  as
the correct generalizations of the $k$-numerical ranges in the present
situation.   We will show below in Theorem 2.6 that in infinite dimensions the
boundary of the numerical range is obtained via a derivative process involving
these sets.  The proof of this result will use various facts about the map
$\Delta$, which we now present.

\begin{theorem} If $0 < s < t < 1$, then the following statements hold.
\begin{enumerate}
\item  The map $\Delta$ is affine on $I_t$.

\item Each  $\Delta_t = \Delta(I_t)$ is a convex, compact subset
of $\mathbb C$.

\item  We have $\displaystyle   \frac{s}{t}I_t \subset I_s$
and so $\Delta_t\subset \Delta_s$.

\item  $\Delta(B\setminus\{0\})= \range(\Delta)$ is convex.

\end{enumerate}
\end{theorem}
\begin{proof} Conclusion $(1)$ is immediate from the definition of
$\Delta$ and the fact that each $I_t$ is convex.  Conclusion $(2)$
follows from conclusion $(1)$ because  $I_t$ is compact and convex and
$\Delta$ is continuous.

If $\xv = (t, z)$ is in $B\setminus \{0\}$, then the convexity of  $B$
implies that the line segment from $\xv$ to $0$ lies in $B$.
Since $s < t$, the point  $\displaystyle \frac{s}{t}\xv$ lies in $B$ and since
$\displaystyle \frac{s}{t}\xv =  \left(s,\frac{s}{t}z\right)$, this point
lies in $I_s$.   Also,
\[
\Delta\left(\frac{s}{t}\xv\right) =\frac{z}{t}= \Delta(\xv).
\]
and so $\Delta_t \subset \Delta_s$.
This proves conclusion $(3)$. Conclusion  $(4)$ follows from assertions
$(2)$ and $(3)$ and the fact that the union of a family of convex sets which is
totally ordered by inclusion is itself a convex set.
\end{proof}

We next record some trivial observations that will be used below.

\begin{proposition} The following statements hold.
\begin{enumerate}
\item If $\gl$ is a complex number, then $W_{ab}(c)-\gl =W_{ab}(c-\gl 1)$

\item If $c^\prime = c-\gl 1$, then the invertible linear map
\[
(x,z) \mapsto (x,z-\gl x)
\]
transforms $B(c)$ onto $B(c^\prime)$.

\item  If   $\tau(c) = 0$, then  the interior of the chord joining
$(0,0)$ and $(1,0)$ (which are boundary points of $B$) is in the relative
interior of
$B$ and $(t,0)$ is in the relative interior  of each isotrace slice $I_t$.
\end{enumerate}
\end{proposition}

Observe that it follows from parts (2) and (3) of Proposition 2.2 that if
$\tau(c) =
\gl$, then the interior of the chord joining $(0,0)$ and $(1,\gl)$ lies in the
relative interior of $B$ and if $0 < t < 1$, then $(t,\gl t)$ lies in the
relative interior of the isotrace slice $I_t$. Thus it is natural to call the
line determined by the origin and $(1,\gl)$ the {\bf central axis} of $B$.

Let us now turn to the definition of the radial complex slopes of
$B(c)$ at the origin.    This concept is best visualized when $\tau(c) = 0$
so that the central axis of $B$ is just the $x$-axis.  In this case we may
informally describe the radial complex slope  of $B$ in the direction $\ta$ as
follows.

Fix an angle $\theta$ and consider the half-plane starting at the x-axis and
making an angle $\theta$ with the positive $y$-axis in the $yz$-plane.  If we
write
$C_\ta$ for the intersection of this half-plane with $B$, then $C_\ta$  is a
convex set of dimension less than or equal to 2. The ``upper" boundary of
$C_\ta$ is the graph of a function $f_{\theta}$ and we define the
radial complex
slope of $B$ at $(0,0)$ in the direction $\ta$ to be the slope of
$f_\ta$ at $0$
as measured in $C_\ta$.

In the general case when $\tau(c) = \gl \ne 0$, the central axis of $B$ is the
line determined by $(1,\gl)$ and the origin, write $c^\prime = c -\gl 1$
and $B^\prime$ for the spectral scale determined by $c^\prime$.  With this we
define the radial complex slope of $B$ at $(0,0)$ in the direction $\ta$ to be
$\mu^\prime_\ta + \gl$, where $\mu^\prime_\ta$ is the radial complex slope of
$B^\prime$ at $(0,0)$.

The  precise definition is as follows.  Suppose $\tau(c) = 0$. If $c = 0$,
then  $B =\{(t,0): 0 \le t\le 1\}$.   In this case the radial complex slope at
$(0,0)$  in the direction $\ta = 0$ is defined to be 0, which is the complex
slope of this chord.

If $c\ne 0$, and the interior of $B$ is empty, then $B$ must be a
planar set whose intersection with the $x$--axis is the interval
$[0,1]$.  In this case we may multiply by an appropriate scalar $\gl =
e^{i\ta}$  and get that $\gl B$ lies in the $(x,y)$--plane. Since we now
have $\tau(\gl ca) \in \mathbb R$ for all $a\in N^+_1$ it follows that
the imaginary part of $\gl c$ is 0 so that $b = \gl c$ is
self--adjoint.  In this
case for each $0 < t < 1$ the associated  isotrace slice of $B(b)$ is a line
segment in the $(x,y)$-plane of the form $\{(t,s): -r_{t,\pi}\le s \le
r_{t,0}\}= \{(t,s): r_{t,\pi}e^{i\pi}\le s\le r_{t,0}\}$.   Hence, we get that
the corresponding isotrace slice of $B(c)$ is the complex interval with
endpoints $r_{t,0}e^{-i\ta}$ and $r_{t,\pi}e^{i(\pi-\ta)}$.  In order to keep
notation consistent with that to be introduced below, we write $r_{t, -\ta} =
r_{t,0}$ and $r_{t,\pi-\ta} = r_{t,\pi}$.  With this we set
\[
r_{-\ta} = \lim_{t\downarrow 0} \frac{r_{t,-\ta}}{t}\text{ and }
r_{\pi-\ta} = \lim_{t\downarrow 0} \frac{r_{t,\pi-\ta}}{t}
\]
and define the radial complex slopes of $c$ at the origin in the
directions $-\ta$ and $\pi-\ta$ to be  $r_{-\ta}e^{-i\ta}$ and
$r_{\pi-\ta}e^{i(\pi-\ta)}$.

Now suppose $B$ has nonempty interior (and $\tau(c) = 0$) so that 0 is
in the interior of each isotrace slice. In this case for each $0 < t
<1$ and
$0 \le \ta < 2\pi$, there is a unique point of the form
$r_{t,\ta}e^{i\ta}$ on the relative boundary of $I_t$.  We define the radial
complex slope of $B$ at $(0,0)$ in the direction $\ta$  to be
$r_\ta e^{i\ta}$, where
\[
r_\ta = \lim_{t\downarrow 0}\frac{r_{t,\ta}}{t}.
\]
We shall show below that each of the limits  above exists and is finite.

Now suppose $c = c_0 + \gl 1$, where $\tau(c_0) = 0$.  In this case we
get that the relative boundary points of each $I_t$ have the form
$r_{t,\ta}e^{i\ta} + \gl t$ by part$(2)$ of Proposition 2.2 and we
define the radial complex slope of $B$ at $(0,0)$ in the direction
$\ta$ to be $r_\ta e^{\it\ta} + \gl$.

Our next goal is to show that the following statements hold.

\begin{enumerate}
\item The boundary points of $W(c)$ are in one to one correspondence with the
complex radial slopes of $B$ at the origin.

\item A point $\gl$ lies in $\partial W(c)\cap W(c)$ if and only if there is
a line segment with complex slope $\gl$ on the boundary of $B$ which is
anchored at the origin.

\item The line segments in $\partial W(c)\cap W(c)$ are in
one to one correspondence with the faces of dimension two on the boundary of
$B$ that contain the origin.
\end{enumerate}
This will be accomplished in several steps. We begin by investigating the
relation between  $\pw{}(c)$ and the geometry of the spectral scale.  It is
useful to begin by presenting two Lemmas.

\begin{lemma}If $\gl$ is in  $\partial W(c)\cap W(c)$, then there is $t >
0$ such that $(t,t\gl) \in B$.
\end{lemma}
\begin{proof} Fix $\gl\in \partial W(c)\cap W(c)$.  Normalizing as in 1.4, we
may assume that $\gl = 0$ and $W(c)$ lies in the right half plane.  Further, we
may define $F_W, r$ and $q$  as in 1.4, except that now we use $b_1$ and $b_2$
instead of $d_1$ and $d_2$.  For example, we now write $r$ for the projection
onto $\script N(b_1)$.

   If $F_W = \{0\}$ or $0$ is an endpoint of $F_W$, then $q > 0$ by
parts (2) and (3) of Lemma 1.5 and we have $\tau(q) = t > 0$ and $cq = 0$
so that $\Psi(q) = (\tau(q),\tau(cq)) = (t,0)$.  Hence, $(t,0) \in B$.

If $F_W$ is one dimensional  and $0$ is in the interior of
$F_W$, then $(rb_2r)^\pm \ne 0$ by part (4) of Lemma 1.5 so that if we write
$r^\pm$ for the range projection of $(rb_2r)^\pm$, then we have $0 < r^\pm < r$
and $b_1 r^\pm = 0$.  Also, since $r^\pm \ne 0$ we get $\tau(r^\pm) > 0$,
$\tau(b_2r^+) >0$ and $\tau(b_2r^-) < 0$.  Thus if we select  suitable small
positive choices of $s^-$ and $s^+$ and write $a = s^--r^- + s^+r^+$, we have
$a \le 1$, $\tau(a) = t > 0$ and $\tau(ca)=0$.  Hence, $\Psi(a) =
(\tau(a),\tau(ca)) = (t,0)$ and so $(t,0) \in B$. Thus in all cases, we get
that there is a point of the form $(t,\gl t) \in B$ with $t > 0$.

\end{proof}

\begin{lemma}If $\gl$ is in  $\partial W(c)$ and there is $t> 0$ such that
$(t,t\gl) \in B$,  then $\gl$ is in $W(c)$.
\end{lemma}

\begin{proof} Fix $\gl\in \partial W(c)$ and assume that there is $t>
0$ such that $(t,t\gl) \in B$ .  Normalizing as in 1.4, we may assume that $\gl
= 0$ and $W(c)$ lies in the right half plane.  Further, we may define $F_W, r$
and $q$ as in 1.4, except that we now use $b_i$'s in place of the $d_i$'s.

With this,  we get that there  is an element $d$ in $N_1^+$ such that
\[
\Psi(d) = (\tau(d), \tau(cd)) =  (\tau(d), \tau(b_1d)+i\tau(b_2d)) = (t,0).
\]
Thus, $\tau(b_1d) = 0$.  Since $W(c)$ lies in the right half plane,
$W(b_1)$ is
nonnegative and so $b_1 \ge 0$. Since $\tau$ is faithful, we get that
$\sqrt{d}b_1\sqrt{d} = 0$ so that $b_1d = 0$.

If $\sqrt{d}b_2\sqrt{d}$ is semi-definite, then since
$\tau(\sqrt{d}b_2\sqrt{d})
=0$, we get $\sqrt{d}b_2\sqrt{d} = 0$ because $\tau$ is faithful.  In this
case since
$d
\ne 0$, we may select a unit vector $\xv$ in the range of $\sqrt{d}$ and get
\[
\lan c\xv,\xv\ran = \lan b_1\xv,\xv\ran + i\lan b_2\xv,\xv\ran = 0
\]
so that $\gl \in W(c)$ in this case.

Finally suppose that $\sqrt{d}b_2\sqrt{d}$ is not semi-definite.  In this case,
we may select unit vectors of the form $\sqrt{d}\xv$ and
$\sqrt{d}\yv$ such that
$\lan b_2 \sqrt{d}\xv,\sqrt{d}\xv\ran =
\lan\sqrt{d}b_2\sqrt{d}\xv,\xv) > 0$ and
$\lan b_2 \sqrt{d}\yv,\sqrt{d}\yv\ran =\lan\sqrt{d}b_2\sqrt{d}\yv,\yv)< 0$.
Since $b_1\sqrt{d} = 0$ and $W(c)$ is convex, we get $\gl = 0 $ in $W(c)$.
\end{proof}

\begin{theorem} We have
\[
\range(\Delta) =  W(c).
\]
\end{theorem}
\begin{proof} We first show that $\overline{\range(\Delta)} = W_{ab}(c)$.
Fix $a\ne 0$ in $N_1^+$ and define the functional $f_a$ on $N$ by the formula
\[
f_a (b) = \tau(ab)/\tau(a).
\]
           Since $\tau$ is a tracial state, $f_a$ is a state of  $N$.  Now
observe that
\[
\Delta(\Psi(a)) = \Delta((\tau(a),\tau(ac)) = \frac{\tau(ac)}{\tau(a)} =
f_a(c).
\]
Thus, $\Delta(\Psi(a)) \in W_{ab}(c)$ and so $\range(\Delta) \subset
W_{ab}(c)$.

To prove that $W_{ab}(c) \subset \overline{\range(\Delta)}$,
we will show that the range of $\Delta$ is dense in $W_{ab}(c)$.  Since both
sets are convex, this will follow if we can show that the closure of the range
of $\Delta$ contains the extreme points of $W_{ab}(c)$. If $\gl$ is any such
extreme point, then by linearity there must be an extreme point  $f$  of  the
set of states on $N$ such that $f(c) = \gl$. Therefore we need only show that
$f$ lies in the weak* closure of the set states of the form $f_a$.  This is an
easy consequence of \cite[Lemma 3.4.1]{Dix}. Hence, $\overline{\range(\Delta)}
= W_{ab}(c)$.

Next, if $\gl \in \partial W(c)\cap W(c)$, then there is $t > 0$ such that
$\xv = (t,t\gl) \in B$ by Lemma 2.3 and so  $\Delta(\xv) = \gl$.  Thus
$\range(\Delta)$ contains $\partial W(c)\cap W(c)$.

If $\gl \in  \partial\range(\Delta)$, then $\gl \in \partial
W(c)$  by the first part of the proof. If in addition $\gl \in \range(\Delta)$,
then there must be a point $(t,t\gl) \in B$ with $t > 0$ and so $\gl \in W(c)$
by Lemma 2.4. Hence
\[
\partial W(c)\cap W(c) = \partial \range(\Delta)\cap\range(\Delta).
\]
Since $W(c)$ and $\range(\Delta)$ are each convex and $\overline{W(c)}=
\overline{\range(\Delta)}$ these sets have the same interiors.  Since they also
have the same relative boundaries, they are equal.

\end{proof}

\begin{theorem}The relative boundary points of $W_{ab}(c)$ are
precisely the radial
complex slopes of $B(c)$ at the origin.
\end{theorem}

\begin{proof} Write $c = c_0 +\gl 1$, where $\tau(c_0) =  0$.  Since
$W_{ab}(c) = W_{ab}(c_0)) + \gl$, by part $(1)$ of Proposition 2.2,
it suffices to establish the Theorem for $c_0$, or in what amounts to
the same thing, we may assume that $\tau(c) = 0$ so that is the $x$-axis is the
central axis of $B$.

If  $c = 0$, then $B = \{(t,0): 0\le t\le 1\}$,  and $W_{ab}(c) =  W(c) = 0$,
which is the radial complex slope of $B$ at $(0,0)$ in the direction  $0$.
Hence, the assertion is true in this case.

Next, if $B$ is a planar set, then we may replace $c$ by a multiple
of itself and assume that $B$ lies in the $(x,y)$-plane. In this case
$c = b$ is
self--adjoint (and $\tau(b) = 0$) and we get that $W_{ab}(b) = [\ga^-
,\ga^+]$, where $\ga^\pm$ are the largest and smallest elements of the spectrum
$b$.  As shown in  \cite[Theorem 1.7]{AAW} there are exactly two slopes of
$B(b)$ at the origin, and these slopes are precisely the largest and smallest
points of the spectrum.  On the other hand, it is clear that the limits defined
above also converge to these slopes and so the Theorem  holds in this case
and it now trivial to extend this to the case where $c = e^{i\ta} b$ and $b$
is self--adjoint.

Now assume that $B$ has nonempty interior, fix $\ta$ and write $r_\ta e^{i\ta}$
for the unique boundary point of $W_{ab}(c)$ with argument $\ta$.   Next
consider the points of the form $(t,r_{t,\ta}e^{i\ta})$ on the boundary of $B$.
Since each $\Delta_t$ is a subset of $W(c)$ by Theorem 2.5, we have
\[
\frac{I_t}{t} \subset W(c)
\]
and therefore $\displaystyle \frac{r_{t,\ta}}{t}\le r_\ta$.
Further, it follows from part $(3)$ of Theorem 2.1 that if $0< s < t <
1$, then
\[
\frac{r_{t,\ta}}{t} <\frac{r_{s,\ta}}{s}.
\]
Since these points lie in $W(c)$ they are bounded by $r_\ta$ so that
\[
r_\ta = \lim_{t\downarrow 0} \frac{r_{t,\ta}}{t}
\]
and so the assertion holds.

\end{proof}

\begin{theorem} A point  $\gl$ is in $\partial W(c) \cap W(c)$ if and only if
there is a line segment on the boundary of $B$ containing the origin and with
complex slope $\gl$.  This segment has the form
\[
\seg_{B,\gl}=\{(t,t\gl): 0 \le t \le t_\gl\},
\]
where $t_\gl > 0$.
\end{theorem}

\begin{proof}    Translating $c$ if necessary, we may  assume that $0$ is in
the relative interior of $W(c)$.  Now  fix $\gl \in \partial W(c) \cap W(c)$ so
that we have $(t_0,\gl t_0) \in B$ for some $t_0 > 0$ by Lemma 2.3 and write
$t_\gl = \sup\{t: (t,t\gl)\in B\}$. Also Write $\gl  = re^{i\ta}$.  If we had
$(t,t\gl)\in B^\circ$ for some $0 < t \le t_\gl$, then we could find $ s > r$
such that if $\mu = se^{i\ta}$, then  $(t,t\mu) \in B$. But in this case we
would have $\mu \in W(c)$by Lemma 2.4  and since $\gl$ is a convex combination
of $\mu$ and $0$, it could not lie on $\partial W(c)$. Hence no such $\mu$ can
exist and so the points $(t,t\gl)$ lies on the boundary of $B$ for each $ 0
\le t  \le t_\gl$. Since a translation preserves the boundary of the spectral
scale, the second assertion is true in the general case.

\end{proof}

In the Theorem below and the sequel we shall be concerned with faces in
$W_{ab}(c), W(c)$ and$B$. In order to avoid confusion we will use $F_{ab}, F_W$
and $F_B$ to denote faces in these sets, respectively.

\begin{theorem} If $F_{ab}$ is a proper one dimensional face of $W_{ab}(c)$ and
we write $F_W = F_{ab}\cap W(c)$, then $F_W$ has dimension one if and only if
there is a two dimensional face $F_B$ in $B$ containing the origin and  such
that each line segment containing the origin in $F_B$ has the form
\[
\seg_{B,\gl} = \{(t,t\gl): 0 \le t \le t_\gl, \gl \in F_W\}.
\]

\end{theorem}

\begin{proof}Suppose that the boundary of $W_{ab}(c)$ contains a face $F_{ab}$
of dimension one  and  $F_W = F_{ab}\cap W(c)$ is one dimensional. Let $\gl_0$
and
$\gl_1$ denote the endpoints of $F_W$ and fix an interior point of the form
$\gl_s = s\gl_0 + (1-s)\gl_1$ for some $0 <s <1$.  Applying  Theorem 2.7, we
get that there is a maximal line segment $\sg_{B,\gl_s}$in the boundary of $B$
containing the origin.  Let $(t_{\gl_s},t_{\gl_s}\gl_s)$ denote the nonzero
endpoint of this line segment so that points of the form $(t, t\gl_s)$  lie on
this line segment for $0 \le t \le t_{\gl_s}$.

Note that an endpoint $\gl_i$ of $F_W$ lies in $F_W$ if an only if there is a
line segment in the boundary of $B$ with complex slope $\gl_i$ and  containing
the origin.  But if we  fix $0 < s_0 < s_1 < 1$,  then  $\sg_{B,\gl_{s_0}}$ and
$\sg_{B,\gl_{s_1}}$ are line segments in the boundary of  $B$. Further if $s_0
\le s \le s_1$  and $t$ is any real number, then
\begin{align*}
\gl_s &= s\gl_0 + (1-s)\gl_1 = \left(\frac{s_1-s}{s_1-s_0}\right)\gl_{s_0} +
\left(1 -\frac{s_1-s}{s_1-s_0}\right)\gl_{s_1}\text{ and } \\
(t,t\gl_s) &= \left(\frac{s_1-s}{s_1-s_0}\right) (t,t\gl_{s_0}) + \left(1 -
\frac{s_1-s}{s_1-s_0}\right)(t,t\gl_{s_1}).\\
\end{align*}
Thus for each $s_0 < s< s_1$, $\seg_{B,s}$ lies in the plane determined by
$\seg_{B,s_0}$ and $\seg_{B,s_1}$  Hence, if we write $t_{\min} =
\min\{t_{\gl_{s_0}}, t_{\gl_{s_1}}\}$, then
\[
P=\{s(t, t\gl_{s_0}) + (1-s)(t,t\gl_{s_1}): s_0 \le  s \le s_1,\, 0 \le t \le
t_{\min}\}
\]
lies on the boundary of $B$.

Since $P$ is a planar set, it must lie in a face $F_B$ of dimension two on the
boundary of $B$ and this face must  contain $\{\sg_{B,\gl_s}: s_0 \le s \le
s_1\}$.  Further, if $0 < s_0 < s_1 < s_2 < 1$, then the same argument shows
that $\{\sg_{B,\gl_s}: s_0 \le s\le s_2\}$ lies in  a two dimensional face
$F_B^\prime$ that contains the planar set
\[
P'=\{s(t, t\gl_0) + (1-s)(t,t\gl_1): s_0 \le  s \le s_2,\, 0 \le t \le
\min(t_{\gl_{s_0}},t_{\gl_{s_2}})\}.
\]
Since $F_B$ and $F_B^\prime$ meet in a planar set, we must have $F_B =
F_B^\prime$.  Hence $F_B$ contains $\{\sg_{B,\gl_s}: 0 < s < 1\}$ and  the
endpoint $\gl_0$ (resp., $\gl_1$) lies in $F_W$, if and only if
$\sg_{B,\gl_0}$ (resp., $\sg_{B,\gl_1}$) lies in $F_B$.
\end{proof}

The following Theorem summarizes the results obtained in Theorems  2.6, 2.7
and 2.8.

\begin{theorem} \label{numranchar} The following statements hold.

\begin{enumerate}
\item The relative boundary points of $W_{ab}(c)$ are precisely the radial
complex slopes of $B(c)$ at the origin.

\item  If $\lambda$ is a relative boundary point of $W_{ab}(c)$, then $\gl \in
W(c)$ if and only if  there is a line segment on the boundary of $B$ that
contains the origin and has  complex slope $\lambda$.

\item If $W_{ab}(c)$ is two dimensional, $F_{ab}$ is a face in $W_{ab}(c)$ of
dimension one and we write $F_W = F_{ab}\cap W(c)$, then  $F_W$ has dimension
one if and only if there is a two dimensional face $F_B$ in $B$ that contains
the origin and such that the complex slopes of the line segments in
$F_B$ that contain the origin consist precisely of the points in $F_W$.

\end{enumerate}
\end{theorem}

\smallskip

As we noted in the introduction to this section, the isotrace slice $I_{1/n}$
is an affine image of $W(c)$ when $N$ acts on a Hilbert space of dimension $n$.
We now show that an analogous result holds in infinite dimensions if the
nonzero extreme points of $B$ are bounded away from $0$.

\begin{corollary}  If
\[
t = \inf\{\tau(p): \Psi(p) \text{ is an extreme point of $B$ and } \Psi(p) \ne
0\} >0,
\]
then the restriction of $\Delta$ to $I_t$ is an affine  isomorphism from $I_t$
onto $W(c)$.
\end{corollary}
\begin{proof} If $t> 0$, so that $B$ has no extreme points of the
form $(s,z)$ with  $0 < s < t$, then the  boundary of $B$ near $0$
consists of line segments joining $0$ and boundary points of $I_t$.

           Since the map $\Delta$ is constant on such line segments, if $0 < s
< t$, then
\[
I_s = \frac{s}{t}I_t.
\]
by part (3) of Theorem 2.1 .  Hence, $\Delta(B\setminus \{0\}) = \Delta_t =
W(c)$ by part $(3)$ and of Theorem 2.1 and Theorem 2.5.
\end{proof}

If $c$ is normal so that $N$ is abelian, then  it was shown in \cite[Theorem
3.4]{AAW} that the spectral scale completely determines $c$ up to
isomorphism.   Nevertheless even in this special case the spectrum of $c$ need
not be visible as a set of slopes (however defined) on the boundary of the
spectral scale.  The content of the next result is that in the normal case the
eigenvalues of $c$ {\em are} visible as the complex slopes of the faces of $B$
of dimension one.  Further, although $\tau$ is not unique in this case,  this
theorem is independent of the choice of $\tau$.

Since the proof of this result uses on the notation and theory developed in
\cite{AAW} and \cite{Geom II}, we now review this material.  If we return to
the view that $B$ is determined by the self-adjoint operators $b_1$ and $b_2$,
then  it was proved in \cite[Theorem
2.3]{AAW} (and revisited in \cite[Corollary 0.3]{Geom II}) that exposed faces
of the spectral scale have the form
\[
F = \Psi([\ptsm,\ptsp]),
\]
where $\tv = (t_1,t_2)$ is a nonzero vector in \rn{2}, $s$ is a real number
and $\ptspm$ are the spectral projections of   $\bt{} = t_1b_1+t_2b_2$
corresponding to the intervals $(-\infty,s]$ and $(-\infty,s)$.

The complete facial structure of a spectral scale was determined in
\cite[Section 3]{Geom II}.  The key to this analysis is the fact that
if $F= \Psi([\ptsm,\ptsp])$ is an exposed face of $B$, then it is an affine
image of a new spectral scale, which is determined as follows.  If we write
$r = \ptsp-\ptsm$ and
\[
A(\xv) = \frac{1}{\tau(r)}(\xv - \Psi(\ptsm)),
\]
then $B_F = A(F)$ is the spectral scale of the cut-down operators $rb_1r$ and
$rb_2r$ restricted to the range of $r$.  This new spectral scale is defined
using the trace $\tau_r$ and the map $\Psi_r$, where,
\[
\tau_r(\cdot) = \frac{1}{\tau(r)}\tau(\cdot)\text{ and } \Psi_r(\cdot) =
\frac{1}{\tau(r)}\Psi(\cdot).
\]

   Thus, if
$F_1$ is an exposed face of $B_F$ of the form $\Psi_r([\qtsmm{1},\qtspp{1}])$,
and
$G$ is the image of
$F_1$ under $A^{-1}$, then
\[
G = \Psi(\ptsm) + \Psi([\qtsmm{1},\qtspp{1}]).
\]

\begin{theorem} If $N$ is abelian, then the point spectrum of $c$ is
exactly the set of complex slopes of the 1-dimensional faces of $B$.
\end{theorem}
\begin{proof} If $F$ is  a 1-dimensional face of $B$, then by \cite [Theorem
3.6]{Geom II}, there are projections $q^- < q^+$ in $N$ such that
$F=\Psi([q^-,q^+])$ and the endpoints of $F$ are  $\Psi(q^\pm)$.   Thus,
the complex slope of $F$ is
\[
\frac{\tau(cq^+)- \tau(cq^-)}{\tau(q^+)- \tau(q^-)}=\frac{\tau(cr)}{\tau(r)}
= \gl.
\]
  Since $F$
has dimension one, if we write $r = q^+-q^-$ then by \cite[Corollary
3.3(1)]{Geom II} $rNr$ has dimension one and therefore $rcr = \gl r$ for some
complex scalar $\gl$. Since $N$ is abelian $rcr = rc = \gl c$ and so $\gl$ is
an eigenvalue for $c$.

Now suppose that $\gl$ is an eigenvalue of $c$.  By using the substitution
$c \rightarrow c-\lambda 1$ we can assume that this eigenvalue is  0.    We
now use the notation and results of \cite{AAW} and \cite{Geom II} as
described prior to the statement of the Theorem.  If we write
$\tv = (1,0)$,  $\bt{} = b_1$ and let $p_{\tv,0}^\pm$ denote the
corresponding spectral projections, then by \cite[Corollary 0.3(1)]{Geom II} $F
=\Psi([p_{\tv,0}^-, p_{\tv,0}^+])$ is an exposed face of $B$. Also,  if we
write  $r = p_{\tv,0}^+ - p_{\tv,0}^-$, then $r \ne 0$ and
$r\bt{} = rb_1 = 0$ by \cite[parts (1) and (3) of Corollary 0.3]{Geom II}.

Now let $p$ denote the projection onto the null space of $b_2$ and observe that
$q = pr$ is the projection onto the null space of $c$, which is nonzero because
$0$ is in the point spectrum of $c$. If $q = r$, then $rNr = rN = qN = qNq
=\mathbb C$ and so $F$ has dimension one by \cite[part (1) of Corollary
3.3]{Geom II}.  Further, the slope of $F$ is
\[
\frac{\tau(\bt{}p_{\tv,0}^+) -\tau(\bt{}p_{\tv,0}^-) }{\tau(p_{\tv,0}^+) -
\tau(p_{\tv,0}^-)} = \frac{\tau(\bt{}r)}{\tau(r)}   = 0,
\]
by \cite[part (3) of Corollary 0.3]{Geom II}.  Hence, the Theorem is true in
this case.

If $0 < q < r$, so that $rNr \ne \mathbb C$, then  $F$ has dimension greater
that one and since $F$ is a proper face, it must have dimension two.  If $B_F =
B(rb_1r,rb_2r)$ as described in the paragraph just before the statement of the
Theorem, we set $\uv = (0,1)$, $s_1 = 0$, and write  $q_{\uv,0}^\pm$
for the corresponding projections, then $\Psi_r([q_{\uv,0}^-, q_{\uv,0}^+]) =
F_1$ is an exposed face of $B_F$.  This face has dimension one  because
$q_{\uv,0}^+ - q_{\uv,0}^-$ is the (nonzero) projection onto the null space of
$rb_2r$.   The corresponding face in $B$ has the form
\[
G = \Psi(p_{\tv,0}^-) +  \Psi([q_{\uv,0}^-,q_{\uv,0}^+]).
\]
Finally, the slope of $G$ is
\[
\frac{\tau(b_2(p_{\tv,0}^- + q_{\uv,0}^+)) - \tau(b_2(p_{\tv,0}^- +
q_{\uv,0}^-))}{\tau(p_{\tv,0}^- + q_{\uv,0}^+) - \tau(p_{\tv,0}^- +
q_{\uv,0}^-)} = \frac{\tau(b_2 q_{\uv,0}^+) -\tau(b_2q_{\uv,0}^-))}{\tau(
q_{\uv,0}^+) - \tau(q_{\uv,0}^-)} = \frac{\tau(cq)}{\tau(q)} = 0.
\]

\end{proof}

\bigskip

\begin{corollary} If $c$ is unitary, then $\sigma(c)$ is exactly the set of
radial complex slopes of
$B$ at $(0,0)$ that have absolute value 1.
\end{corollary}
\begin{proof} If $\gl  \in \gs(c)$, then $|\gl| = 1$ and so it lies in the
boundary of the numerical range.  Hence  $\gl$ is a radial complex slope of $B$
by Theorem 2.6.

If $\gl$ is a radial complex slope of $B$ at $(0,0)$ with $|\gl| =1$,
then $\gl$
lies in the boundary of $W_{ab}(c)$ by Theorem 2.6.  Since $|\gl| = 1$ it is an
extreme point of this set and so $\gl \in \gs(c)$ by Theorem 1.7.
\end{proof}

\end{document}